\documentclass[12pt,twoside,a4paper]{article}
 
\usepackage{bm}
\usepackage{amsmath,amssymb,amsthm}
\usepackage{enumerate}
\usepackage{fullpage}
\usepackage{pslatex} 
\usepackage{verbatim}
\usepackage{url}

\newcommand{\Odip}[2]{\mathcal{O}_{#1}\!\Bigl(#2\Bigr)\mathchoice{\!}{}{}{}}

\newcommand{\odip}[2]{{o}_{#1}\!\left(#2\right)\mathchoice{\!}{}{}{}}
\newcommand{\odi}[1]{\odip{}{#1}}

\newcommand{\Eulerphi}{\varphi}
\newcommand{\dx}{\mathrm{d}}

\newcommand{\NN}{\mathfrak{N}}

\newcommand{\M}{\mathfrak{M}}
\newcommand{\m}{\mathfrak{m}}
\renewcommand{\t}{\mathfrak{t}}
\newcommand{\Z}{\mathbb{Z}}
\newcommand{\Q}{\mathbb{Q}}
\newcommand{\R}{\mathbb{R}}
\newcommand{\N}{\mathbb{N}}
\newcommand{\singseries}{\mathfrak{S}}

\newtheorem{Lemma}{Lemma}

\newenvironment{Proof}[1][Proof]{\par\noindent\textbf{#1.}~}
  {\hfill$\square$\smallskip\par}

  \newtheoremstyle{Nonumtheorems}
  {10pt}
  {6pt}
  {\itshape}
  {}
  {\bfseries}
  {.}
  {.5em}
  {\thmname{#1}\thmnote{ (#3)}}

\theoremstyle{Nonumtheorems}
\newtheorem{Nonumthm}{Theorem}
\newtheorem{Nonumcor}{Corollary}

\begin{document}

\title{On a Diophantine problem with two primes and $s$ powers of two}
\author{A.~LANGUASCO and A.~ZACCAGNINI}
 
\date{}
\maketitle

\begin{abstract}
We refine a recent result of Parsell \cite{Parsell2003}
on the values of the form
$
\lambda_1p_1
+
\lambda_2p_2
+
\mu_1 2^{m_1}
+
\dotsm
+
\mu_s 2^{m_s},
$
where $p_1,p_2$ are prime numbers, $m_1,\dotsc, m_s$ are positive
integers, $\lambda_1 / \lambda_2$ is negative and irrational and
$\lambda_1 / \mu_1$, $\lambda_2/\mu_2 \in \Q$. \par
\medskip
\noindent
2000 AMS Classification: 11D75, 11J25, 11P32, 11P55. \par
\end{abstract}

\section{Introduction}
\allowdisplaybreaks

In this paper we are  interested to study the values of the form
\begin{equation}
\label{linear-form}
\lambda_1p_1
+
\lambda_2p_2
+
\mu_1 2^{m_1}
+
\dotsm
+
\mu_s 2^{m_s},
\end{equation}
where $p_1,p_2$ are prime numbers, $m_1,\dotsc, m_s$
are positive integers, and the coefficients $\lambda_1$, $\lambda_2$ and
$\mu_1$, \dots, $\mu_s$ are real numbers satisfying suitable relations.

This is clearly a variation of the so-called Goldbach-Linnik problem,
\emph{i.e.} to prove that every sufficiently large even integer is
a sum of two primes and $s$ powers of two,
where $s$ is a fixed integer. 
Concerning this problem the first result was proved by Linnik himself
\cite{Linnik51,Linnik53} who remarked that a suitable $s$ exists
but he gave no explicitly estimate of its size. Other results
were proved by Gallagher \cite{Gallagher1975},
Liu, Liu \& Wang \cite{LiuLW98a,LiuLW98b, LiuLW99},
Wang \cite{Wang99} and Li \cite{Li2000, Li2001}.
Now the best conditional result is due to Pintz \& Ruzsa
\cite{PintzR2003} and Heath-Brown \& Puchta \cite{Heath-BrownP2002}
($s=7$ suffices under the assumption of the Generalized Riemann Hypothesis),
while, unconditionally, it is due to Heath-Brown \& Puchta
\cite{Heath-BrownP2002} ($s = 13$ suffices).
Elsholtz, in unpublished work, improved it to $s=12$.
We should also remark that Pintz \& Ruzsa announced a proof for the
case $s=8$ in their paper \cite{PintzR2006} which is as yet
unpublished.
Looking for the size of the exceptional set of the Goldbach problem
we recall the fundamental paper
by Montgomery-Vaughan \cite{MontgomeryVaughan1975}
in which they showed that the number of even integers up to $X$ that
are not the sum of two primes is $\ll X^{1-\delta}$.
Pintz recently announced that $\delta=1/3$ is admissible
in the previous estimate. Concerning the
exceptional set for the  Goldbach-Linnik problem, the authors
of this paper in a joint work with Pintz \cite{LanguascoPZ2007}
proved that for every $s \ge 1$, there
are $\ll X^{3/5} (\log X)^{10}$ even integers in $[1,X]$ that are not
the sum of two primes and $s$ powers of two.
This obviously corresponds to the case
$\lambda_1 = \lambda_2 = \mu_1 = \dots = \mu_s = 1$.

In diophantine approximation several results were proved
concerning the linear forms with primes that, in some
sense, can be considered as the real analogous
of the binary and ternary Goldbach problems.
On this topic we recall the papers by
Vaughan \cite{Vaughan1974a,Vaughan1974b,Vaughan1976},
Harman \cite{Harman1991},
Br\"udern-Cook-Perelli \cite{BrudernCP1997},
and Cook-Harman \cite{CookH2006}.

Concerning the problem in \eqref{linear-form}, we can consider it
as a real analogous of the Goldbach-Linnik problem.
We have the following
\begin{Nonumthm}
Suppose that $\lambda_1,\lambda_2$ are real numbers
such that $\lambda_1/\lambda_2$ is negative and irrational
with $\lambda_1>1$, $\lambda_2<-1$ and
$\vert \lambda_1/\lambda_2 \vert \geq1$.
Further suppose that $\mu_1, \dotsc,  \mu_s$ are nonzero
real numbers such that 
$\lambda_i/\mu_i \in \Q$, for $i\in\{1,2\}$, and denote by
$a_i/q_i$  their reduced representations as rational numbers.
Let moreover $\eta$ be a sufficiently small positive constant such that
$\eta<\min(\lambda_1/a_1;\vert \lambda_2/a_2 \vert)$.
Finally,
for $\lambda_1/\lambda_2$ transcendental, let
\begin{equation}
\label{s0-def-transc}
s_0
=
2
+
\Bigl\lceil
\frac
{
\log (2C(q_1,q_2)\vert \lambda_1 \lambda_2\vert)
-
\log \eta
}
{-\log (0.91237810306)}
\Bigr\rceil,
\end{equation}
while, for $\lambda_1/\lambda_2$ algebraic, let
\begin{equation}
\label{s0-def-alg}
s_0
=
2+
\Bigl\lceil
\frac
{
\log (2C(q_1,q_2)\vert \lambda_1 \lambda_2\vert)
-
\log  \eta
}
{-\log (0.83372131685)}
\Bigr\rceil,
\end{equation}
where $C(q_1,q_2)$ 
verifies
\begin{equation}
\label{Cq1q2-def}  
C(q_1,q_2) 
=
\Bigl(
\log{2}
+
C \cdot  
\singseries^{\prime}(q_1) 
\Bigr)^{1/2}
\Bigl(
\log{2}
+
C \cdot  
\singseries^{\prime}(q_2) 
\Bigr)^{1/2}, 
\end{equation}
with
\begin{equation}
\label{singseries-def}
  \singseries^{\prime}(n)
  =
  \prod_{\substack{p \mid n \\ p > 2}} \frac{p - 1}{p - 2}
\end{equation}
and
$C =  10.0219168340$.

\noindent
Then  for every real number $\gamma$ and every integer
$s\geq s_0$ the  inequality
\begin{equation}
\label{main-inequality}
\vert
\
\lambda_1p_1
+
\lambda_2p_2
+
\mu_1 2^{m_1}
+
\dotsm
+
\mu_s 2^{m_s}
+
\gamma
\
\vert
<
\eta
\end{equation}
has infinitely many solutions in primes
 $p_1,p_2$ and positive integers $m_1,\dotsc, m_s$.
\end{Nonumthm}
The only result on this problem we know is by Parsell \cite{Parsell2003}; 
our values in \eqref{s0-def-transc}-\eqref{s0-def-alg} 
improve Parsell's 
one given by
\begin{equation}
\label{Parsell-1}
s_0
=
2+
\Bigl\lceil
\frac
{
\log (2C_1(q_1,q_2)\vert \lambda_1 \lambda_2\vert)
-
\log  \eta
}
{-\log (0.954)}
\Bigr\rceil,
\end{equation} 
where 
\begin{equation}
\label{Parsell-2}
C_1(q_1,q_2)
=
25 (\log 2q_1)^{1/2}(\log 2q_2)^{1/2}.
\end{equation}
Checking the proof in \cite{Parsell2003} one can see that \eqref{Parsell-2}
is in fact 
\begin{equation}
\label{Parsell-3}
C_1(q_1,q_2,\epsilon)
=
\Bigl(
1+ C_1  \cdot 
\singseries^{\prime}(q_1) 
\Bigr)^{1/2}
\Bigl(
1+ C_1  \cdot 
\singseries^{\prime}(q_2)
\Bigr)^{1/2}
+  \epsilon,
\end{equation}
and $C_1  = 11.4525218267$.
Comparing the numerical values involved in 
\eqref{s0-def-transc}-\eqref{Cq1q2-def} with \eqref{Parsell-1} and \eqref{Parsell-3},
without considering the contribution of the $\log 2$
which in \eqref{Parsell-3} is replaced by $1$,
we see that the our gain is 
about  $50$\%
in the transcendental case and about  $75$\%
in the algebraic case.
For instance, taking $\lambda_1= \sqrt{3}= \mu_1^{-1}$, 
$\lambda_2= -\sqrt{2}= \mu_2^{-1}$ and $\eta=1$, we get $s_0= 61$
while  for $\lambda_1= \pi = \mu_1^{-1}$, 
$\lambda_2= -\sqrt{2}= \mu_2^{-1}$ and $\eta=1$, we get $s_0= 119$.
In both cases,  \eqref{Parsell-1} gives $s_0=267$.

Moreover we remark that the  work of Rosser-Schoenfeld 
\cite{RosserS1962} on $n/\Eulerphi(n)$, see Lemma \ref{sing-series-estim} below,
gives for $\singseries'(q)$ a sharper estimate than $2\log (2q)$, 
used in \eqref{Parsell-2}, for large values of $q$.

With respect to \cite{Parsell2003}, 
our main gain comes 
from enlarging the size of the major arc since this lets us
use sharper estimates on the minor arc.
In particular, on the major arc we replaced the technique used in 
\cite{Parsell2003} 
with a well-known argument involving the
Selberg integral; this also simplified the actual
work to get a ``good'' major arc contribution.

On the minor arc we used Br\"udern-Cook-Perelli's \cite{BrudernCP1997}
and Cook-Harman's \cite{CookH2006} technique to deal with the exponential sum 
on primes ($S(\alpha)$)  while, in order to work with the exponential sum over powers of two
($G(\alpha)$),  we inserted Pintz-Ruzsa's \cite{PintzR2003} algorithm to estimate 
the measure of the subset of the minor arc on which $\vert G(\alpha)\vert$ is ``large''.
These two ingredients lead to
a sharper estimate on the minor arc 
and let us improve the size of the denominators in 
\eqref{s0-def-transc}-\eqref{s0-def-alg}.
It is in this step that we have to distinguish whether  $\lambda_1/\lambda_2$ 
is an algebraic or a transcendental number; this fact leads to
two different estimates for the minor arc and, \emph{a fortiori},
using Pintz-Ruzsa's algorithm (see Lemma \ref{minor-arc-power-of-two-estim}), 
to two different constants in \eqref{G-algebraic}-\eqref{G-transcendental}
and \eqref{s0-def-transc}-\eqref{s0-def-alg}.

A second, less important, gain arises from our Lemma \ref{Dioph-equation}
which improves the values
in \eqref{Cq1q2-def} comparing with
the ones in  \eqref{Parsell-3} (obtained in \cite{Parsell2003}, Lemma 3).
Such an improvement 
comes from using  the Prime Number Theorem 
(to get $\log 2$ instead of $1$) 
and Khalfalah-Pintz's \cite{KhalfalahP2006}
computational estimates for the number of representations
of an integer as a difference of powers of two, see Lemma \ref{KP-Lemma}.

Finally we remark that assuming a suitable form of the twin-prime conjecture, 
\emph{i.e.} $B=1$ in Lemma \ref{BD-Thm2}, we get that \eqref{Cq1q2-def}
holds with $C= 2.5585042082$.

Using the notation
$\bm{\lambda}=(\lambda_1,\lambda_2)$,
$\bm{\mu}=(\mu_1,\mu_2)$, as a consequence of the Theorem we have the 
\begin{Nonumcor} 
Suppose that $\lambda_1,\lambda_2$ are real numbers
such that $\lambda_1/\lambda_2$ is negative and irrational.
Further suppose $\mu_1, \dotsc,  \mu_s$ are nonzero
real numbers such that $\lambda_i/\mu_i \in \Q$, for $i\in\{1,2\}$, 
and denote by $a_i/q_i$  their reduced representations as rational numbers.
Let moreover $\eta$ be a sufficiently small positive constant such that
$\eta<\min(\vert \lambda_1/a_1\vert ;\vert \lambda_2/a_2 \vert)$ and 
$\tau\geq\eta>0$.
Finally let $s_0=s_0(\bm{\lambda},\bm{\mu},\eta)$ as defined 
in \eqref{s0-def-transc}-\eqref{s0-def-alg}.
Then
for every real number $\gamma$ and every integer
$s\geq s_0$ the  inequality
\begin{equation}
\label{general-inequality}
\vert
\
\lambda_1p_1
+
\lambda_2p_2
+
\mu_1 2^{m_1}
+
\dotsm
+
\mu_s 2^{m_s}
+
\gamma
\
\vert
<
\tau
\end{equation}
has infinitely many solutions in primes
 $p_1,p_2$ and positive integers $m_1,\dotsc, m_s$.
\end{Nonumcor}

This Corollary immediately follows from the Theorem 
since, multiplying by a suitable constant both sides
of \eqref{general-inequality}, we can always reduce ourselves
to study  the case $\lambda_1>1$, $\lambda_2<-1$ and
$\vert \lambda_1/\lambda_2 \vert \geq1$.
Hence the Theorem assures us that \eqref{main-inequality}
has infinitely many solutions and the Corollary immediately 
follows from the condition $\tau\geq\eta$.

We finally remark that the condition 
about about the rationality of the two ratios $\lambda_i/\mu_i$,
$i=1,2$,
which, at first sight, could appear a ``weird'' one,
is in fact quite natural in the sense that otherwise the numbers $\lambda x+
\mu y$, $x,y\in \Z$, are dense in $\R$ by Kronecker's Theorem, see
also the remark after Lemma \ref{Dioph-equation}.

\paragraph{Acknowledgements.}
We would like to thank J\'{a}nos Pintz, Umberto Zannier and Carlo Viola
for a discussion, Imre Ruzsa
for sending us his original U-Basic code for
Lemma \ref{minor-arc-power-of-two-estim}
and Karim Belabas for helping us in improving the performance
of our PARI/GP code for the Pintz-Ruzsa algorithm.

\section{Definition}

Let $\epsilon$ be a sufficiently small positive constant, $X$ be a large parameter, 
$M=\vert \mu_1 \vert + \dotsm +\vert \mu_s \vert$
and $L=\log_2 (\epsilon X/(2M))$, where $\log_2 v$ is
the base $2$ logarithm of $v$.
We will use the Davenport-Heilbronn  variation of the Hardy-Littlewood
method to count the number of solutions
$\NN(X)$ of the inequality \eqref{main-inequality}  with
$\epsilon X \leq p_1, p_2 \leq X$
and $1\leq m_1, \dotsc, m_s \leq L$.
Let now $e(u) = \exp(2\pi i u)$ and
\[
S(\alpha) = \sum_{\epsilon X\leq p \leq X} \log p\ e(p\alpha)
\quad
\textrm{and}
\quad
G(\alpha)= \sum_{1 \leq m \leq L}  e(2^m\alpha).
\]
For $\alpha\neq 0$, we also define
\begin{equation}
\notag
K(\alpha,\eta)=
\Bigl(\frac{\sin\pi \eta \alpha}{\pi \alpha}\Bigr)^2
\end{equation}
and hence both 
\begin{equation}
\label{K-transform}
\widehat{K}(t,\eta)
=
\int_\R K(\alpha,\eta) e(t\alpha ) \dx\alpha
=
\max(0; \eta -\vert t\vert)
\end{equation}
and 
\begin{equation}
\label{K-inequality}
K(\alpha,\eta)
\ll
\min(\eta^2; \alpha^{-2})
\end{equation}
are well-known facts. 
Letting
\begin{equation}
\notag
I(X ; \R)
=
\int_{\R}
S(\lambda_1 \alpha) S(\lambda_2 \alpha)
G(\mu_1 \alpha)\dotsm G(\mu_s \alpha)
e(\gamma \alpha)
K(\alpha,\eta)
\dx \alpha,
\end{equation}
it follows from \eqref{K-transform} that
\[
I(X ; \R)
  \ll
\eta \log^2 X \cdot \NN(X).
\]
We will prove that
\begin{equation}
\label{I-lower-bound}
I(X ; \R)
\gg_{s,\bm{\lambda},\epsilon}
\eta^2 X (\log X)^s
\end{equation}
thus obtaining
\[
 \NN(X) \gg_{s,\bm{\lambda},\epsilon} \eta  X (\log X)^{s-2}
\]
and hence the Theorem follows. 
To prove \eqref{I-lower-bound}
we first dissect the real line in the major, minor and trivial arcs,
by choosing  $P=X^{1/3}$ and letting
\begin{equation}
\label{dissect-def}
\M = \{\alpha\in \R: \vert \alpha \vert \leq P/X \}, 
\quad 
\m = \{\alpha\in \R: P/X<\vert \alpha \vert \leq L^2\},
\end{equation}
and $\t = \R \setminus(\M\cup\m)$.
Accordingly, we  write
\begin{equation}
\label{integral-dissect}
I(X ; \R)
=
I(X ; \M) + I(X ; \m) + I(X ; \t).
\end{equation}
We will prove that the inequalities
\begin{equation}
\label{major-goal}
I(X ; \M)
\geq
c_1  \eta^2 X L^s,
\end{equation}
\begin{equation}
\label{trivial-goal}
\vert 
I(X ; \t)
\vert
=
\odi{XL^s},
\end{equation}
hold for all sufficiently large $X$, and
\begin{equation}
\label{minor-goal}
\vert
I(X ; \m)
\vert
\leq
c_2(s) \eta X L^s,
\end{equation}
where $c_2(s)>0$ depends on $s$,
$c_2(s)\to 0$ as $s\to+\infty$,  and
$c_1=c_1(\epsilon,\bm{\lambda})>0$ is a constant such that
\begin{equation}
\label{constant-condition}
c_1 \eta -c_2(s) \geq c_3\eta
\end{equation}
for some absolute positive constant  $c_3$ and $s\geq s_0$.
Inserting \eqref{major-goal}-\eqref{constant-condition} into \eqref{integral-dissect},
we finally obtain that \eqref{I-lower-bound} holds thus proving the Theorem.

\section{Lemmas}

Let $1\leq n \leq (1-\epsilon)X/2$ be an integer and $p,p^{\prime}$ two prime numbers.
We define the twin prime counting function as follows
\begin{equation}
\label{Z-def}
Z(X; 2n)
=
\sum_{\epsilon X \leq p \leq X}
\sum_{\substack{p^{\prime}\leq X \\ p^{\prime} -p =2n}} \log p \log p^{\prime}.
\end{equation}
Moreover we denote by $\singseries(n)$ the singular series and set
$\singseries(n) = 2 c_0 \singseries^{\prime}(n)$ where $\singseries^{\prime}(n)$
is defined in \eqref{singseries-def} and 
\begin{equation}
\label{c0-def}
   c_0
  =
   \prod_{p > 2} \Big( 1 - \frac{1}{(p-1)^2} \Big).
\end{equation}
Notice that $\singseries^{\prime}(n)$ is a multiplicative function.
According to Gourdon-Sebah \cite{GourdonS2001},
we can also write that
$0.66016181584<c_0<0.66016181585$.

Let further $k\geq 1$ be an integer and $r_{k,k}(m)$ be the
number of representations of an integer $m$ as
$\sum_{i=1}^{k} 2^{u_i} - \sum_{i=1}^{k} 2^{v_i}$, where
$1\leq u_i, v_i \leq L$ are integers, so that
$r_{k,k}(m)=0$ for sufficiently large $\vert m\vert$.
Define
\[
S(k,L) =
\sum_{m\in \Z\setminus \{0\}}
r_{k,k}(m)
\singseries(m).
\]
The first Lemma is about the behaviour of $S(k,L)$
for sufficiently large  $X$.
\begin{Lemma}[Khalfalah-Pintz \cite{KhalfalahP2006}, Theorem 2]
\label{KP-Lemma}
For any given $k\geq 1$,
there exists $A(k)\in \R$ such that
\[
\lim_{L\to +\infty}
\Bigl(
\frac{S(k,L)}{2L^{2k}}
-1
\Bigr)
=
A(k).
\]
\end{Lemma}
Moreover they also proved numerical estimates
for $A(k)$ when $1\leq k\leq 7$. We will just need
\begin{equation}
\label{A(1)-estim}
A(1) <  0.2792521041.
\end{equation}
The second lemma is an upper bound for the multiplicative part of
the singular series.
\begin{Lemma}
\label{sing-series-estim} 
For $n\in \N$, $n\geq 3$, we have that
\[
\singseries^{\prime}(n)
<
\frac{n}{c_0\Eulerphi(n)}
<
\frac{e^{\gamma} \log \log n}{c_0}
+
\frac{2.50637}{c_0 \cdot \log \log n},
\]
where $\gamma=0.5772156649\dotsc$ is the Euler constant.
\end{Lemma}
\begin{Proof} 
Let $n\geq 3$. 
The first estimate follows immediately remarking
that
\[
\singseries^{\prime}(n)
=
\prod_{\substack{p \mid n \\ p > 2}}
\frac{(p - 1)^2}{p(p - 2)}
\prod_{\substack{p \mid n \\ p > 2}}
\frac{p}{p - 1}
<
\prod_{\substack{p > 2}}
\frac{(p - 1)^2}{p(p - 2)}
\prod_{\substack{p \mid n}}
\frac{p}{p - 1}
=
\frac{1}{c_0}
\frac{n}{\Eulerphi(n)}.
\]
The second estimate is a direct application of Theorem 15 of
Rosser and Schoenfeld \cite{RosserS1962}.
\end{Proof}
Letting $f(1)=f(2)=1$ and 
$ f(n) = n/(c_0\Eulerphi(n)) $
for  $n \geq 3$,
we can say that the inequality $\singseries^{\prime}(n) \leq f(n)$
is sharper than Parsell's estimate  $\singseries^{\prime}(n) \leq 2 \log(2 n)$,
see page 7 of \cite{Parsell2003},  for every $n\geq 1$.
Since it is clear that computing the exact value of $f(n)$ for
large values of $n$ it is not easy (it requires the knowledge of
every prime factor of $n$), we also remark that
the second estimate in Lemma \ref{sing-series-estim} leads to
a sharper bound  than $\singseries^{\prime}(n) \leq 2 \log(2 n)$ 
for every $n\geq 14$.

The next lemma is a famous result
of Bombieri and Davenport.
\begin{Lemma}[Theorem 2 of Bombieri-Davenport \cite{BombieriD1966}]
\label{BD-Thm2}
There exists a positive constant $B$ such that,
for every positive integer $n$, we have
\[
Z(X; 2n)
<
B\  \singseries(n) X,
\]
where $Z(X; 2n)$ and $\singseries(n)$ are defined in
\eqref{singseries-def} and \eqref{Z-def}-\eqref{c0-def}, 
provided that $X$ is sufficiently large.
\end{Lemma}
Chen \cite{Chen1978} proved that $B =3.9171$
can be used in Lemma \ref{BD-Thm2}. 
The assumption of a suitable form of the twin prime conjecture,
\emph{i.e.} $Z(X; 2n) \sim \singseries(n) X$
for $X\to +\infty$,  implies that in this case we can take $B=1$.

Now we state some lemmas we need to estimate  $I(X ;\m)$.
The first one is
\begin{Lemma}
\label{Dioph-equation}
Let $X$ be a sufficiently large parameter and let $\lambda, \mu \neq
0$ be two real numbers such that $\lambda / \mu\in \Q$.
Let $a,q\in \Z\setminus\{0\}$ with $q>0$, $(a,q)=1$ be such
that $\lambda/\mu= a/q$.
Let further $0<\eta < \vert \lambda/a \vert $.
We have
\[
\int_{\R}
\vert
S(\lambda \alpha) G(\mu \alpha)
\vert^2
K(\alpha, \eta)
\dx\alpha
<
\eta X L^2
\Bigl(
(1-\epsilon)\log{2}
+
C \cdot 
\singseries'(q) 
\Bigr)
+
\Odip{M,\epsilon}{\eta X L},
\]
where $C = 10.0219168340$. 
\end{Lemma}
\begin{Proof}
First of all we remark that the constant 
$C $ is in fact $2 B  (1+A(1))$,
where  $B = 3.9171$ is the constant
in Lemma \ref{BD-Thm2} and $A(1)$ is
estimated in \eqref{A(1)-estim}.
This
should be compared with the value
$C_1 =11.4525218267$  
obtained in \cite{Parsell2003}. 
Assuming $B=1$ in Lemma \ref{BD-Thm2}, we get
$C=2.5585042082$.
Letting now
\[
I
=
\int_{\R}
\vert
S(\lambda \alpha) G(\mu \alpha)
\vert^2
K(\alpha, \eta)
\dx\alpha,
\]
by \eqref{K-transform} we immediately have
\begin{equation}
\label{expanded}
I
=
\sum_{\epsilon X  \leq p_1, p_2\leq X}
\sum_{1 \leq m_1, m_2 \leq L}
\log p_1 \log p_2
\max
\Bigl(
0;
\eta - \vert
\lambda(p_1-p_2) +\mu(2^{m_1}-2^{m_2})
\vert
\Bigr).
\end{equation}
Let $\delta= \lambda(p_1-p_2) + \mu(2^{m_1}-2^{m_2})$.
For a sufficiently small $\eta>0$, we claim that
\begin{equation}
\label{delta0}
\vert
\delta
\vert
< \eta
\quad
\textrm{is equivalent to}
\quad
\delta = 0.
\end{equation}
Recall our hypothesis on $a$ and $q$, and assume that
$\delta\neq 0$ in \eqref{delta0}.
For $\eta<\vert \lambda/a\vert$ this leads to a contradiction. In fact we have
\[
\frac{1}{\vert a \vert}
>
\frac{\eta}{\vert \lambda\vert}
>
\Bigl\vert
(p_1-p_2) + \frac{q}{a}(2^{m_1}-2^{m_2})
\Bigr\vert
=
\Bigl\vert
\frac{a(p_1-p_2) + q(2^{m_1}-2^{m_2})}{a}
\Bigr\vert
\geq
\frac{1}{\vert a \vert},
\]
since $a(p_1-p_2) + q(2^{m_1}-2^{m_2})\neq 0$ is a linear integral
combination.
Inserting \eqref{delta0} in \eqref{expanded},
for $\eta<\vert \lambda/a\vert$ we can write that
\begin{equation}
\label{expanded1}
I
=
\eta
\sum_{\epsilon X  \leq p_1, p_2\leq X}
\sum _{\substack{
1\leq m_1, m_2 \leq L \\ \hskip-1.7cm
\lambda(p_1-p_2) +\mu(2^{m_1}-2^{m_2}) =0
}}
\log p_1 \log p_2.
\end{equation}
The diagonal contribution  in \eqref{expanded1}
is equal to
\begin{equation}
\label{diagonal-contrib}
\eta
\sum_{\epsilon X \leq p \leq X}
\log^2 p
\sum_{1\leq m \leq L}
1
=
\eta
XL^2
(1-\epsilon)\log{2}
+\Odip{M,\epsilon}{\eta X L}
\end{equation}
where we used the Prime Number Theorem instead
of trivially estimate the contribution of $\log p_i$
as in \cite{Parsell2003}.

Now we have to estimate
the contribution $I'$ of the non-diagonal solutions
of $\delta=0$ and we will achieve this by connecting $I'$ with
the singular series of the twin prime problem.
Recalling that $\lambda/\mu = a/q \neq 0$, $(a,q)=1$,
by Lemma \ref{BD-Thm2} and the fact that
$Z(X;(q/a)(2^{m_2}-2^{m_1}))\neq0$
if and only if
$a \mid (2^{m_2}-2^{m_1})$, we have,
since $\singseries(v)=\singseries(2^{u}v)$
for every $u,v\in \N$, $u\geq1$, that
\begin{equation}
\label{Sol-estim1}
\begin{split}
I' 
\leq
2  \eta
\sum_{1\leq m_1 < m_2 \leq L}
Z
\Bigl(
X;\frac{q}{a}(2^{m_2}-2^{m_1})
\Bigr) 
<
2 
B 
X  \eta
\sum_{1\leq m_1 < m_2 \leq L}
\singseries
\Bigl(
\frac{q}{a}(2^{m_2}-2^{m_1})
\Bigr).
\end{split}
\end{equation}
Using the multiplicativity of
$\singseries^{\prime} (n)$ (defined in \eqref{singseries-def}),
we get
\[
\singseries^{\prime}
\Bigl(
\frac{q}{a}(2^{m_2}-2^{m_1})
\Bigr)
\leq
\singseries^{\prime}(q)
\singseries^{\prime}
\Bigl(
\frac{2^{m_2}-2^{m_1}}{a}
\Bigr)
\leq
\singseries^{\prime}(q)
\singseries^{\prime}(2^{m_2}-2^{m_1})
\]
and so, by Lemma \ref{KP-Lemma},
\eqref{A(1)-estim} and \eqref{Sol-estim1},
we can write, for every sufficiently large $X$, that
\begin{equation}
\label{I-estim2}
\begin{split}
I'
&
\leq
2 
B 
X  \eta \singseries^{\prime}(q)
\sum_{1\leq m_1 < m_2 \leq L}
\singseries(2^{m_2}-2^{m_1})
=
B 
X  \eta
\singseries^{\prime}(q)
S(1,L) 
\\
&
<
2B (1+A(1))
\singseries^{\prime}(q)
X  \eta L^2.
\end{split}
\end{equation}
Hence, by \eqref{expanded1}-\eqref{diagonal-contrib}
and \eqref{I-estim2}, we finally get
\[
I
<
\eta X L^2
\Bigl(
(1-\epsilon)\log{2}
+
2B (1+A(1))
\singseries'(q) 
\Bigr)
+
\Odip{M,\epsilon}{\eta X L},
\]
this way proving
Lemma \ref{Dioph-equation}.
\end{Proof}

We remark that if in Lemma \ref{Dioph-equation} we consider also the case
$\lambda/\mu \in \R\setminus\Q$,
we can just find $\eta=\eta(X)\to 0$ as  $X\to+\infty$ and this
implies that $s_0\approx |\log \eta| \to +\infty$,
see equations \eqref{s0-def-transc}-\eqref{s0-def-alg}
for the precise definition of $s_0$. This essentially depends on the fact
that, for $\lambda/\mu\in\R\setminus\Q$ and $m,n\in \Z$, it
is not possible to find a function $f(X)$ such that
$\vert \lambda m  + \mu n \vert
\geq f(X)$ and $f(X)\to c>0$ as $X \to +\infty$ since
the set of values of $\lambda m + \mu n$
is dense in $\R$.
A different, but related, way to see this phenomenon is to remark
that the inequality $|\alpha n + m| < \eta$ is equivalent to the
pair of inequalities $\Vert n \alpha \Vert < \eta$ or
$\Vert n \alpha \Vert > 1 - \eta$, where $\Vert u \Vert$
is the distance of $u$ from the nearest integer.
When $\alpha$ is irrational, it has $\sim 2 \eta X$ solutions with
$n \le X$, since the sequence $\Vert n \alpha \Vert$ is uniformly
distributed modulo 1.

To estimate the contribution of $G(\alpha)$ on the minor arc
we use Pintz-Ruzsa's method as developed in
\cite{PintzR2003}, \S 3-7.
\begin{Lemma}[Pintz-Ruzsa \cite{PintzR2003}, \S~7]
\label{minor-arc-power-of-two-estim}
Let $0< c <1$. Then there exists $\nu=\nu(c)\in (0,1)$ 
such that 
\[
\vert 
E(\nu) 
\vert
:=
\vert 
\{
\alpha \in (0,1) \
\textrm{such that} \ 
\vert 
G(\alpha) 
\vert 
> \nu L
\}
\vert
\ll_{M,\epsilon}
X^{-c}.
\]
\end{Lemma}
To obtain explicit values for $\nu$ we had to write 
our own version of Pintz-Ruzsa algorithm since
in this application the estimates has to be performed
for a different choice of parameters than the ones
they used in \cite{PintzR2003}.
We used the PARI/GP \cite{PARI2} scripting language and the
gp2c compiling tool to be able to compute fifty decimal digits
(but we write here just ten) of the constant
involved in the following Lemma.
We will write two different estimates that we will use in the case
$\lambda_1/\lambda_2$ is a transcendental or an algebraic
number. 
Running the program in our cases, Lemma \ref{minor-arc-power-of-two-estim}
gives the following results:
\begin{equation}
\label{G-algebraic}
\vert
G(\alpha)
\vert
\leq
0.83372131685 \cdot L
\end{equation}
if
$\alpha \in [0,1]  \setminus E$ where
$\vert E \vert \ll_{M,\epsilon} X^{-2/3-10^{-20}}$, to be
used when $\lambda_1 / \lambda_2$ is algebraic,
and 
\begin{equation}
\label{G-transcendental}
\vert
G(\alpha)
\vert
\leq
0.91237810306   \cdot L
\end{equation}
if
$\alpha \in [0,1] \setminus E$ where
$\vert E \vert \ll_{M,\epsilon} X^{-4/5-10^{-20}}$, to be
used when $\lambda_1 / \lambda_2$ is transcendental.

The computing time to get \eqref{G-algebraic}-\eqref{G-transcendental} 
on a double quad-core PC of the NumLab laboratory of the
Department of Pure and Applied Mathematics of the 
University of Padova  was equal in the first case to
24 minutes and 40 seconds (but to get 30 correct digits just 
3 minutes and 24 seconds suffice) and to 29 minutes
(but to get 30 correct digits just 3 minutes and 50 seconds suffice)
in the second case.
You can download the PARI/GP source code of our program
together with the cited numerical values at the following link:
\url{www.math.unipd.it/~languasc/PintzRuzsaMethod.html}.

Now we state some lemmas we will use to work on the major arc.
Let $\theta(x)=\sum_{p\leq x} \log p$,
\begin{equation}
\label{Selberg-int-def}
J(X,h) 
=
\int_{\epsilon X}^X
(\theta(x+h)- \theta(x) -h)^2
\dx x
\end{equation}
be the Selberg integral and
\[
  U(\alpha)
  =
    \sum_{\epsilon X\leq n \leq X} e(\alpha n).
\]

Applying a famous Gallagher's lemma (\cite{Gallagher1970}, Lemma 1)
on the truncated $L^2$-norm
of exponential sums to $S(\alpha) - U(\alpha)$, one gets the following
well-known statement which we cite from 
Br\"udern-Cook-Perelli \cite{BrudernCP1997}, Lemma 1.
\begin{Lemma}
\label{BCP-Gallagher}
For $1/X \leq Y \leq 1/2$ we have
\[
\int_{-Y}^Y
\vert
S(\alpha) - U(\alpha)
\vert^2
\dx \alpha
\ll_{\epsilon}
\frac{\log X}{Y}
+
Y^2X
+
Y^2 J \Bigl( X,\frac{1}{Y} \Bigr),
\]
where $J(X,h)$ is defined in \eqref{Selberg-int-def}.
\end{Lemma}
To estimate the Selberg integral, we use the next result.
\begin{Lemma}[Saffari-Vaughan \cite{SaffariV1977a}, \S6]
\label{Saffari-Vaughan}
For any $A>0$ there exists $B=B(A)>0$ such that
\[
J(X,h)
\ll_{\epsilon}
\frac{h^2X}
{(\log X)^A}
\]
uniformly for $h\geq X^{1/6} (\log X)^B$.
\end{Lemma}

\section{The major arc}
Letting
\begin{equation}
\label{T-def-estim}
  T(\alpha)
  =
 \int_{\epsilon X}^{X}e(t\alpha)\dx t
  \ll_{\epsilon}
  \min \Bigl(X, \frac{1}{\vert\alpha\vert} \Bigr),
\end{equation}
we first  write
\begin{equation}
\label{I-splitting}
\begin{split}
I(X ; \M)
&=
\int_\M
T(\lambda_1 \alpha)
T(\lambda_2 \alpha)
G(\mu_1\alpha)
\dotsm
G(\mu_s\alpha)
e(\gamma \alpha)
K(\alpha, \eta) \dx \alpha
\\
& +
\int_\M
\Bigl(S(\lambda_1 \alpha) - T(\lambda_1 \alpha)\Bigr)
T(\lambda_2 \alpha)
 G(\mu_1\alpha)
\dotsm
G(\mu_s\alpha)
e(\gamma \alpha)
K(\alpha, \eta) \dx \alpha \\
& +
\int_\M
S(\lambda_1 \alpha)
\Bigl(S(\lambda_2 \alpha) - T(\lambda_2 \alpha)\Bigr)
G(\mu_1\alpha)
\dotsm
G(\mu_s\alpha)
e(\gamma \alpha)
K(\alpha, \eta) \dx \alpha \\
&=
J_1+ J_2+J_3,
\end{split}
\end{equation}
say.
In what follows we will prove that
\begin{equation}
\label{J1-lower-bound}
J_1
\geq
\frac{1-(7/2)\lambda_1\epsilon}
{2 \vert\lambda_1\lambda_2 \vert}
\eta^2 XL^s
\end{equation}
and
\begin{equation}
\label{J2-estim}
J_2 + J_3
=
\odi{\eta^2XL^s},
\end{equation}
thus obtaining by \eqref{I-splitting}-\eqref{J2-estim} that
\[
I(X ; \M)
\geq
\frac{1-4\lambda_1\epsilon}
{2 \vert\lambda_1\lambda_2 \vert}
\eta^2 XL^s.
\]
Thus we will prove that \eqref{major-goal} holds with
$c_1= (1-4\lambda_1\epsilon)/ (2 \vert\lambda_1\lambda_2 \vert)$.

\paragraph{Estimation of $J_2$ and $J_3$.}
We first estimate $J_3$.
We remark that, by the partial summation formula, we have
$T(\alpha) - U(\alpha) \ll (1 + X\vert \alpha \vert)$.
So, recalling $P=X^{1/3}$, \eqref{dissect-def} 
and $\vert S(\lambda_1\alpha) \vert \ll X \log X$,
we get
\[
\int_{\M}
\vert
T(\lambda_2\alpha)
-
U(\lambda_2\alpha)
\vert
\vert S(\lambda_1\alpha) \vert
\dx \alpha
 \ll
X \log X
\int_{\M}
(1 + X\vert \lambda_2 \alpha \vert)
\dx \alpha
\ll_{\bm{\lambda}}
X^{2/3} \log X.
\]

Hence,
using the  trivial estimates
$\vert G(\mu_i \alpha)\vert \leq L$, $K(\alpha,\eta)\ll \eta^2$, we can write
\[
J_3
=
\int_\M
S(\lambda_1 \alpha)
\Bigl(S(\lambda_2 \alpha) - U(\lambda_2 \alpha)\Bigr)
G(\mu_1\alpha)
\dotsm
G(\mu_s\alpha)
e(\gamma \alpha)
K(\alpha, \eta) \dx \alpha 
+
\Odip{\bm{\lambda},M}
{\eta^2X^{2/3} L^{s+1}}.
\]

Now using \eqref{dissect-def}, the Cauchy-Schwarz inequality, the Prime Number Theorem, 
Lemmas \ref{BCP-Gallagher}-\ref{Saffari-Vaughan}
with $A=3$, $Y=P/X$, $P=X^{1/3}$, and again
the  trivial estimates
$\vert G(\mu_i \alpha)\vert \leq L$, $K(\alpha,\eta)\ll \eta^2$,
we have that
\[
\begin{split}
J_3
&
\ll
\eta^2 L^s
\Bigl(
\int_{\M}
\vert
S(\lambda_2\alpha)
-
U(\lambda_2\alpha)
\vert^2
\dx \alpha
\Bigr)^{1/2}
\Bigl(
\int_{\M}
\vert
S(\lambda_1\alpha)
\vert^2
\dx \alpha
\Bigr)^{1/2}
+
\Odip{\bm{\lambda}, M}
{\eta^2X^{2/3} L^{s+1}}
\\
&
\ll_{\bm{\lambda}, M, \epsilon}
\eta^2 L^s
\frac{X^{1/2}}
{(\log X)^{3/2}}
\Bigl(
\int_{0}^{1}
\vert
S(\alpha)
\vert^2
\dx \alpha
\Bigr)^{1/2}
+
\eta^2X^{2/3} L^{s+1}
\ll_{\bm{\lambda},M,\epsilon}
\eta^2 X L^{s-1}
=
\odi{\eta^2XL^s}.
\end{split}
\]

The integral $J_2$ can be estimated analogously using \eqref{T-def-estim}
instead of the Prime Number Theorem.
Hence  \eqref{J2-estim} holds.

\paragraph{Estimation of $J_1$.}
Recalling that $P= X^{1/3}$
and using  \eqref{dissect-def}, \eqref{T-def-estim} and \eqref{I-splitting}
we obtain
\begin{equation}
\label{J1J-relation}
J_1
=
\sum_{1\leq m_1 \leq L} \dotsm \sum_{1\leq m_s \leq L}
J
\Bigl(
\mu_1 2^{m_1} + \dotsm + \mu_s 2^{m_s}
+ \gamma
\Bigr)
+
\Odip{\epsilon}{\eta^2  X^{2/3} L^s},
\end{equation}
where $J(u)$ is defined by
\[
J(u)
:=
\int_\R
T(\lambda_1 \alpha)
T(\lambda_2 \alpha)
e(u \alpha)
K(\alpha, \eta) \dx \alpha
=
  \int_{\epsilon X}^{X}
  \int_{\epsilon X}^{X}
\widehat{K}(\lambda_1u_1+\lambda_2u_2 + u)
\dx u_1 \dx u_2
\]
and the second relation follows by interchanging the order of integration.
Assume
now that $\vert u \vert \leq \epsilon X$ and that
$2\epsilon\lambda_1 X
\leq
\vert \lambda_2\vert u_2
\leq
(1-\epsilon\lambda_1)X$.
For
$\eta<2\epsilon(\lambda_1-1)X$
and
$X$ sufficiently large, we have, by \eqref{K-transform},
that there exists an interval
for $u_1$, of length $\geq \eta/\lambda_1$ and contained
in $[\epsilon X,X]$, on which
$\widehat{K}(\lambda_1u_1+\lambda_2u_2 + u) \geq \eta/2$.
Thus we have
\begin{equation}
\label{J-lower-bound}
J(u)
\geq
\frac{1-3\lambda_1\epsilon}
{2 \vert\lambda_1\lambda_2 \vert}
\eta^2 X.
\end{equation}
For a sufficiently large $X$, it is clear that
$\vert \mu_1 2^{m_1} + \dotsm + \mu_s 2^{m_s}
+ \gamma
\vert \leq \epsilon X$
while the other condition on the size of
$\vert \lambda_2\vert u_2$ follows from the hypothesis
$\vert \lambda_1/\lambda_2\vert \geq 1$
and
$\lambda_2<-1$.
Hence,
from \eqref{J1J-relation}-\eqref{J-lower-bound},
we obtain that \eqref{J1-lower-bound} holds.

\section{The trivial arc}

Recalling \eqref{dissect-def}, the trivial estimate 
$\vert G(\mu_i \alpha)\vert \leq L$
and using the Cauchy-Schwarz inequality, we get
\[
\vert I(X ; \t) \vert
\ll
L^s
\Bigl(
\int_{L^2}^{+\infty}
\vert
S(\lambda_1\alpha)
\vert^2
K(\alpha,\eta)
\dx \alpha
\Bigr)^{1/2}
\Bigl(
\int_{L^2}^{+\infty}
\vert
S(\lambda_2\alpha)
\vert^2
K(\alpha,\eta)
\dx \alpha
\Bigr)^{1/2}
\]
By \eqref{K-inequality}
and making a change of variable, we have,
for $i=1,2$, that
\[
\begin{split}
\int_{L^2}^{+\infty}
\vert
S(\lambda_i\alpha)
\vert^2
K(\alpha,\eta)
\dx \alpha
& 
\ll_{\bm{\lambda}}
\int_{\lambda_i L^2}^{+\infty}
\frac
{\vert
S(\alpha)
\vert^2
}
{\alpha^2}
\dx \alpha
\ll
\sum_{n\geq \lambda_i L^2}
\frac{1}{(n-1)^2}
\int_{n-1}^{n}
\vert
S(\alpha)
\vert^2
\dx \alpha
\\
&
\ll_{\bm{\lambda}}
L^{-2}
\int_0^{1}
\vert
S(\alpha)
\vert^2
\dx \alpha
\ll_{\bm{\lambda},M,\epsilon}
\frac{X}{\log X},
\end{split}
\]
by the Prime Number Theorem,
and hence \eqref{trivial-goal} holds.

\section{The minor arc: $\lambda_1/\lambda_2$ algebraic}

Recalling first
\[
I(X ; \m)
=
\int_{\m}
S(\lambda_1 \alpha) S(\lambda_2 \alpha)
G(\mu_1 \alpha)\dotsm G(\mu_s \alpha)
e(\gamma \alpha)
K(\alpha,\eta)
\dx \alpha,
\]
and letting $c\in (0,1)$ to be chosen later,
we first split $\m$ as  $\m_1 \cup \m_2$, 
$\m_1 \cap \m_2=\emptyset$, 
where $\m_2$ is the set of $\beta\in \m$
such that $\vert G(\beta)\vert > \nu(c) L$
and $\nu(c)$ is defined in Lemma  \ref{minor-arc-power-of-two-estim}.
We will choose $c$ to get
$\vert I(X ; \m_2) \vert = \odi{\eta X}$,
since, again by  Lemma  \ref{minor-arc-power-of-two-estim},
we know that $\vert \m_2 \vert \ll_{M,\epsilon} s L^2 X^{-c}$.

To this end, we first use the trivial estimates
$\vert G(\mu_i \alpha)\vert \leq L$ and $K(\alpha,\eta)\ll \eta^2$,
and the Cauchy-Schwarz inequality thus obtaining
\begin{equation}
\label{minor-arc-2}
\begin{split}
\vert I(X ; \m_2) \vert
& \leq
L^{s}
\Bigl(
\int_{\m_2}
\vert S(\lambda_1 \alpha) S(\lambda_2 \alpha) \vert ^2
K(\alpha,\eta)
\dx \alpha
\Bigr)^{1/2}
\Bigl(
\int_{\m_2}
K(\alpha,\eta)
\dx \alpha
\Bigr)^{1/2}
\\
&
\ll
\eta L^s
\vert \m_2 \vert ^{1/2}
\Bigl(
\int_{\m_2}
\vert S(\lambda_1 \alpha) S(\lambda_2 \alpha) \vert ^2
K(\alpha,\eta)
\dx \alpha
\Bigr)^{1/2}.
\end{split}
\end{equation}
We can now argue as in section 4 of 
Br\"udern-Cook-Perelli \cite{BrudernCP1997}
thus getting
\begin{equation}
\label{minor-arc-2-low}
\int_{\m_2}
\vert S(\lambda_1 \alpha) S(\lambda_2 \alpha) \vert ^2
K(\alpha,\eta)
\dx \alpha
  \ll_{\epsilon}
\eta X^{8/3+\epsilon'}.
\end{equation}
Hence, by \eqref{minor-arc-2-low}, 
\eqref{minor-arc-2} becomes
\[
\vert I(X ; \m_2) \vert
\ll_{M,\epsilon}
s^{1/2} \
\eta^{3/2}
X^{4/3+2\epsilon'-c/2}.
\]
Taking $c=2/3+10^{-20}$ and using \eqref{G-algebraic}, we get,
for $\nu= 0.83372131685$ and a
sufficiently small $\epsilon'>0$, that
\begin{equation}
\label{minor-arc-2-final}
\vert I(X ; \m_2) \vert
=
\odi{\eta X}.
\end{equation}

Now we evaluate the contribution of $\m_1$.
Using Lemma \ref{Dioph-equation}
and the  Cauchy-Schwarz inequality,
we have
\begin{align} 
\vert I(X ; \m_1) \vert
&
\leq
(\nu L)^{s-2}
\notag
\\
& \times
\Bigl(
\int_{\m}
\vert S(\lambda_1 \alpha) G(\mu_1 \alpha) \vert ^2
K(\alpha,\eta)
\dx \alpha
\Bigr)^{1/2}
\Bigl(
\int_{\m}
\vert S(\lambda_2 \alpha) G(\mu_2 \alpha) \vert ^2
K(\alpha,\eta)
\dx \alpha
\Bigr)^{1/2}
\notag
\\
&
<
\nu^{s-2} 
C(q_1,q_2)
 \eta X L^{s},
\label{minor-arc-1}
\end{align}
where, recalling Lemmas \ref{sing-series-estim} and
\ref{Dioph-equation},
$C(q_1,q_2)$ is defined as we did in \eqref{Cq1q2-def}.

Hence, by \eqref{minor-arc-2-final} and  \eqref{minor-arc-1},
for $X$ sufficiently large
we finally get
\begin{equation} 
\notag
\vert I(X ; \m) \vert
<
(0.83372131685)^{s-2}
C(q_1,q_2)
 \eta X L^{s}
\end{equation}
whenever $\lambda_1/\lambda_2$ is an algebraic number.
This means that \eqref{minor-goal} holds, in this case, with
$c_2(s)=(0.83372131685)^{s-2}
C(q_1,q_2)$.

\section{The minor arc: $\lambda_1/\lambda_2$ transcendental}

We will act on $\m_1$ as in \eqref{minor-arc-1}
of the previous section thus obtaining
\begin{equation}
\label{minor-arc-1-transc}
\vert I(X ; \m_1) \vert
<
\nu^{s-2} 
C(q_1,q_2)
\eta X
L^{s},
\end{equation}
where $C(q_1,q_2)$ is defined in
\eqref{Cq1q2-def}.

Now we proceed to estimate $I(X ; \m_2)$.
First we argue as in the previous section until \eqref{minor-arc-2}
and then we work as in section 8 of Cook-Harman \cite{CookH2006}
and pp.~221-223 of Harman \cite{Harman1991} thus obtaining
\[
\int_{\m_2}
\vert S(\lambda_1 \alpha) S(\lambda_2 \alpha) \vert ^2
K(\alpha,\eta)
\dx \alpha
\ll
\eta^2
X^{14/5+\epsilon'}
+
\eta
X^{13/5+\epsilon'}.
\]
This, using \eqref{minor-arc-2}, leads to 
\[
\vert I(X ; \m_2) \vert
\ll_{M,\epsilon}
s^{1/2}
X^{-c/2}
(
\eta^{2}
X^{7/5+\epsilon'}
+
\eta^{3/2}
X^{13/10+\epsilon'}
).
\]
Taking $c=4/5+10^{-20}$ and using \eqref{G-transcendental}, we get,
for $\nu= 0.91237810306$ and a
sufficiently small $\epsilon'>0$,
that
\begin{equation}
\label{minor-arc-2-transc}
\vert I(X ; \m_2) \vert
=
\odi{\eta X}.
\end{equation}

Hence, by  \eqref{minor-arc-1-transc} and \eqref{minor-arc-2-transc},
for $X$ sufficiently large
we finally get
\begin{equation} 
\notag
\vert I(X ; \m) \vert
<
(0.91237810306)^{s-2}
C(q_1,q_2)
 \eta X L^{s}
\end{equation}
whenever $\lambda_1/\lambda_2$ is a transcendental number.
This means that \eqref{minor-goal} holds, in this case, with
$c_2(s)=(0.91237810306)^{s-2}
C(q_1,q_2)$.

\section{Proof of the Theorem} 

We have to verify if there exists an $s_0\in\N$ such that
\eqref{constant-condition} holds for $X$ sufficiently large.
Combining the inequalities \eqref{major-goal}-\eqref{minor-goal},
where
$c_2(s)= (0.83372131685)^{s-2}C(q_1,q_2)$
if $\lambda_1/\lambda_2$ is algebraic and,
if $\lambda_1/\lambda_2$ is transcendental,
$c_2(s)= (0.91237810306)^{s-2}C(q_1,q_2)$,
we obtain for $s\geq s_0$, $s_0$ defined 
in \eqref{s0-def-transc}-\eqref{s0-def-alg},
that \eqref{constant-condition} holds in both cases.
This completes the proof of the Theorem.

\normalsize
\bigskip
\begin{tabular}{l@{\hskip 14mm}l}
A.~Languasco               & A.~Zaccagnini \\
Universit\`a di Padova     & Universit\`a di Parma \\
Dipartimento di Matematica & Dipartimento di Matematica \\
Pura e Applicata           & Parco Area delle Scienze, 53/a \\
Via Trieste 63             & Campus Universitario \\
35121 Padova, Italy        & 43100 Parma, Italy \\
{\it e-mail}: languasco@math.unipd.it & {\it e-mail}:
alessandro.zaccagnini@unipr.it
\end{tabular}

\end{document}